\documentclass[reqno]{amsproc}
\usepackage{amssymb}
\usepackage{euscript}
\usepackage{mathrsfs}
\usepackage{units}
\usepackage{color}

\makeatletter
\@namedef{subjclassname@2010}{%
\textup{2010} Mathematics Subject Classification}
\makeatother

\numberwithin{equation}{section}

\newtheorem{thm}{Theorem}[section]

\newtheorem{lem}[thm]{Lemma}
\newtheorem{pro}[thm]{Proposition}

\newtheorem*{thm*}{Theorem}

\theoremstyle{remark}
\newtheorem{rem}[thm]{Remark}

\theoremstyle{definition}

\DeclareMathOperator{\lin}{\mbox{\sc lin}}
\DeclareMathOperator{\D}{d}
\DeclareMathOperator{\dess}{{\mathsf{Des}}}
\DeclareMathOperator{\dzii}{{\mathsf{Chi}}}
\DeclareMathOperator{\E}{e}
\DeclareMathOperator{\koo}{{\mathsf{root}}}

\DeclareMathOperator{\paa}{{\mathsf{par}}}

\newcommand*{\cbb}{\mathbb C}

\newcommand*{\des}[1]{{\dess(#1)}}
\newcommand*{\dz}[1]{{\EuScript D}(#1)}
\newcommand*{\dzi}[1]{\dzii(#1)}

\newcommand*{\dzin}[2]{\dzii^{\langle#1\rangle}(#2)}

\newcommand*{\escr}{{\mathscr{E}_V}}
\newcommand*{\ff}{\mathcal F}
\newcommand*{\Ge}{\geqslant}
\newcommand*{\hh}{\mathcal H}

\newcommand*{\I}{{\mathrm i}}
\newcommand*{\is}[2]{\langle#1,#2\rangle}

\newcommand*{\kk}{\mathcal K}

\newcommand*{\lambdab}{{\boldsymbol\lambda}}

\newcommand*{\Le}{\leqslant}

\newcommand*{\nbb}{\mathbb N}
\newcommand*{\ogr}[1]{\boldsymbol B(#1)}
\newcommand*{\ob}[1]{{\mathcal R}(#1)}
\newcommand*{\pa}[1]{\paa(#1)}

\newcommand*{\slam}{S_{\boldsymbol \lambda}}
\newcommand*{\smalloplus}{\raise0pt\hbox{$\scriptscriptstyle \oplus$}}

\newcommand*{\tcal}{{\mathscr T}}
\newcommand*{\tbb}{{\mathbb T}}

\newcommand*{\xx}{\mathcal X}
\newcommand*{\zbb}{\mathbb Z}

\begin{document}
   \title[Normal extensions
escape from the class of weighted shifts] {Normal
extensions escape from the class of weighted shifts \\
on directed trees}
   \author[Z.\ J.\ Jab{\l}o\'nski]{Zenon Jan Jab{\l}o\'nski}
\address{Instytut Matematyki, Uniwersytet Jagiello\'nski,
ul.\ \L ojasiewicza 6, PL-30348 Kra\-k\'ow, Poland}
   \email{Zenon.Jablonski@im.uj.edu.pl}
   \author[I.\ Jung]{Il Bong Jung}
   \address{Department of Mathematics, Kyungpook National
University, Daegu 702-701, Korea}
   \email{ibjung@knu.ac.kr}
   \author[J.\ Stochel]{Jan Stochel}
\address{Instytut Matematyki, Uniwersytet Jagiello\'nski,
ul.\ \L ojasiewicza 6, PL-30348 Kra\-k\'ow, Poland}
   \email{Jan.Stochel@im.uj.edu.pl}
   \thanks{Research of the first
and the third authors was supported by the MNiSzW
(Ministry of Science and Higher Education) grant NN201
546438 (2010-2013). The second author was supported by
Basic Science Research Program through the National
Research Foundation of Korea (NRF) funded by the
Ministry of Education, Science and Technology
(2009-0087565).}
   \dedicatory{Dedicated to Professor Franciszek H.
Szafraniec on the occasion of his 70th birthday}
   \begin{abstract}
A formally normal weighted shift on a directed tree is
shown to be a bounded normal operator. The question of
whether a normal extension of a subnormal weighted
shift on a directed tree can be modeled as a weighted
shift on some, possible different, directed tree is
answered.
   \end{abstract}
   \maketitle
   \section{Introduction}
The notion of a weighted shift on a directed tree has
been introduced and studied extensively in
\cite{j-j-s}. As shown therein, the class of weighted
shifts on directed trees is wide enough to contain
operators with some subtle properties including
hyponormal operators whose squares are not hyponormal
and non-hyponormal paranormal operators. It is
well-known that there are no (classical) unilateral or
bilateral weighted shifts with the aforesaid
properties. In two recent papers \cite{j-j-s3,j-j-s4}
new examples of unbounded operators with pathological
properties have been constructed, each of them being
implemented as a weighted shift on a special directed
tree. The former contains an example of a hyponormal
operator whose square has trivial domain, the latter
an example of a non-hyponormal operator the
$C^\infty$-vectors of which generate Stieltjes moment
sequences.

In \cite{j-j-s,b-j-j-s} the question of subnormality
of bounded and unbounded weighted shifts on directed
trees has been studied. Criteria for subnormality of
such operators written in terms of consistent systems
of probability measures have been established. When
analyzing subnormality, a question arises as to
whether a normal extension of a nonzero subnormal
weighted shift on a directed tree $\tcal$ with nonzero
weights can be modeled (up to unitary equivalence) as
a weighted shift on some, possibly different, directed
tree. As shown in Section \ref{nevws}, in most
instances this is not the case. The only exceptional
cases are those in which the directed tree $\tcal$ is
isomorphic to either $\zbb$ or $\zbb_+$ (cf.\ Theorem
\ref{normwsno}). In the latter case the normal
extension can be modeled as a weighted shift on a
directed tree which comes from $\zbb$ by gluing a leaf
to the directed tree $\zbb$ at the vertex $0$ (cf.\
Remark \ref{rem}).
   \section{Preliminaries} In what follows,
$\zbb$ stands for the set of all integers and $\cbb$
for the set of all complex numbers. We also use the
following notation
   \begin{align*}
\text{$\zbb_+ = \{n \in \zbb\colon n \Ge 0\}$ and
$\nbb = \{n \in \zbb\colon n\Ge 1\}$.}
   \end{align*}
Let $A$ be an operator in a complex Hilbert space
$\hh$ (all operators considered in this paper are
linear). Denote by $\dz{A}$, $\ob{A}$ and $A^*$ the
domain, the range and the adjoint of $A$ (in case it
exists). A closed densely defined operator $N$ in
$\hh$ is said to be {\em normal} if $N^*N=NN^*$
(equivalently:\ $\dz{N}=\dz{N^*}$ and
$\|N^*h\|=\|Nh\|$ for all $h \in \dz{N}$). For this
and other facts concerning unbounded operators we
refer the reader to \cite{b-s,weid}. A densely defined
operator $S$ in $\hh$ is said to be {\em subnormal} if
there exists a complex Hilbert space $\kk$ and a
normal operator $N$ in $\kk$ such that $\hh \subseteq
\kk$ (isometric embedding) and $Sh = Nh$ for all $h
\in \dz S$; such an $N$ is called a {\em normal
extension} of $S$. We refer the reader to \cite{con}
for the theory of bounded subnormal operators,
\cite{s-sz1,s-sz2,s-sz3,StSz2} for the foundations of
the theory of unbounded subnormal operators and
\cite{c-j-k,kou,k-t1,k-t2} for research related to
special classes of subnormal operators. From now on,
$\ogr \hh$ stands for the $C^*$-algebra of all bounded
operators $A$ in $\hh$ such that $\dz{A}=\hh$. We
write $\lin \ff$ for the linear span of a subset $\ff$
of $\hh$.

Let $\tcal=(V,E)$ be a directed tree ($V$ stands for
the set of all vertices of $\tcal$ and $E$ for the set
of all edges of $\tcal$). If $\tcal$ has a root, which
will always be denoted by $\koo$, then we write
$V^\circ:=V\setminus \{\koo\}$; otherwise, we put
$V^\circ = V$. Set $\dzi u = \{v\in V\colon (u,v)\in
E\}$ for $u \in V$. If for a given vertex $u \in V$
there exists a unique vertex $v\in V$ such that
$(v,u)\in E$, then we denote it by $\pa u$. The
correspondence $u \mapsto \pa u$ is a partial function
from $V$ to $V$. For $n \in \nbb$, the $n$-fold
composition of the partial function $\paa$ with itself
will be denoted by $\paa^n$. Let $\paa^0$ stand for
the identity map on $V$. We say that $\tcal$ is {\em
leafless} if $V =\{u \in V \colon \dzi u \neq
\varnothing\}$. A vertex $u \in V$ is called a {\em
branching vertex} of $\tcal$ if $\dzi{u}$ consists of
at least two vertices. For a subset $W$ of $V$, we
define $\dzi W = \bigcup_{v \in W} \dzi v$ and $\des W
= \bigcup_{n=0}^\infty \dzin n W$, where
   \allowdisplaybreaks
    \begin{gather*}
\dzin{0}{W} = W, \quad \dzin{n+1}{W} =
\dzi{\dzin{n}{W}}, \quad n\in \zbb_+.
    \end{gather*}
For $u \in V$, we put $\dzin n u = \dzin n {\{u\}}$
and $\des{u}=\des{\{u\}}$. The functions
$\dzin{n}{\cdot}$ and $\des{\cdot}$ have the following
properties (see e.g., \cite[Proposition
2.2.1]{b-j-j-s}):
   \begin{align}  \label{num4}
\dzin{n}{u} &= \{w \in V\colon \paa^n(w)=u\}, \quad n
\in \zbb_+,\, u \in V,
   \\
\des u & = \bigsqcup_{n=0}^\infty \dzin n u, \quad
u\in V, \label{num3}
   \\
\des{u_1} \cap \des{u_2} & = \varnothing, \quad u_1,
u_2 \in \dzi{u},\, u_1 \neq u_2,\, u \in V,
\label{num3+}
   \end{align}
where the symbol $\bigsqcup$ is reserved to denote
pairwise disjoint union of sets.

   Let $\ell^2(V)$ be the Hilbert space of all square
summable complex functions on $V$ equipped with the
standard inner product. For $u \in V$, we define $e_u
\in \ell^2(V)$ to be the characteristic function of
the one point set $\{u\}$. The family $\{e_u\}_{u\in
V}$ is an orthonormal basis of $\ell^2(V)$; we call it
the {\em canonical orthonormal basis} of $\ell^2(V)$.
Set
   \begin{align*}
\escr = \lin \{e_u\colon u \in V\}.
   \end{align*}
Given $\lambdab = \{\lambda_v\}_{v \in V^\circ}
\subseteq \cbb$, we define the operator $\slam$ in
$\ell^2(V)$ by
   \begin{align*}
   \begin{aligned}
\dz {\slam} & = \{f \in \ell^2(V) \colon
\varLambda_\tcal f \in \ell^2(V)\},
   \\
\slam f & = \varLambda_\tcal f, \quad f \in \dz
{\slam},
   \end{aligned}
   \end{align*}
where $\varLambda_\tcal$ is the map defined on
functions $f\colon V \to \cbb$ via
   \begin{align*}
(\varLambda_\tcal f) (v) =
   \begin{cases}
\lambda_v \cdot f\big(\pa v\big) & \text{ if } v\in
V^\circ,
   \\
0 & \text{ if } v=\koo.
   \end{cases}
   \end{align*}
$\slam$ is called a {\em weighted shift} on the
directed tree $\tcal$ with weights $\{\lambda_v\}_{v
\in V^\circ}$. Note that any weighted shift $\slam$ on
$\tcal$ is a closed operator (cf.\ \cite[Proposition
3.1.2]{j-j-s}). Combining Propositions 3.1.3, 3.1.8,
3.4.1 and 3.1.7 of \cite{j-j-s}, we get the following
properties of $\slam$ (from now on, we adopt the
convention that $\sum_{v\in\varnothing} x_v=0$).
   \begin{pro}\label{bas}
Let $\slam$ be a weighted shift on a directed tree
$\tcal$ with weights $\lambdab = \{\lambda_v\}_{v \in
V^\circ}$. Then the following assertions hold{\em :}
   \begin{enumerate}
   \item[(i)] $e_u$ is in $\dz{\slam}$ if and only if
$\sum_{v\in\dzi u} |\lambda_v|^2 < \infty$; if $e_u
\in \dz{\slam}$, then
   \begin{align} \label{eu}
\slam e_u = \sum_{v\in\dzi u} \lambda_v e_v \quad
\text{and} \quad \|\slam e_u\|^2 = \sum_{v\in\dzi u}
|\lambda_v|^2,
   \end{align}
   \item[(ii)] $\slam$ is densely defined if and only
if $\escr \subseteq \dz{\slam}$,
   \item[(iii)]  $\slam \in \ogr{\ell^2(V)}$ if and only
if $\alpha_{\lambdab}:=\sup_{u\in
V}\sum\nolimits_{v\in\dzi u} |\lambda_v|^2 < \infty$;
moreover, if $\slam \in \ogr{\ell^2(V)}$, then
$\|\slam\|^2=\alpha_{\lambdab}$,
   \item[(iv)] if $\slam$ is densely defined, then
$\escr \subseteq \dz{\slam^*}$ and
   \begin{align} \label{sl*}
\slam^*e_u=
   \begin{cases}
   \overline{\lambda_u} e_{\pa u} & \text{if } u \in V^\circ, \\
0 & \text{if } u = \koo,
   \end{cases}
   \quad u \in V,
   \end{align}
   \item[(v)] $\slam$ is injective if and only if
$\tcal$ is leafless and $\sum_{v\in\dzi u}
|\lambda_v|^2 > 0$ for every $u\in V$.
   \end{enumerate}
   \end{pro}
   \section{Formal normality - a general structure}
Recall that a densely defined operator $N$ in a
complex Hilbert space $\hh$ is said to be {\em
formally normal} if $\dz{N} \subseteq \dz{N^*}$ and
$\|N^*h\|=\|N h\|$ for all $h \in \dz{N}$ (cf.\
\cite{Cod0,Cod}). In this section we show that
formally normal weighted shifts on directed trees are
always bounded and normal.
   \begin{pro}\label{normal}
If $\slam$ is a nonzero weighted shift on a directed
tree $\tcal$ with weights $\lambdab = \{\lambda_v\}_{v
\in V^\circ}$, then the following three conditions are
equivalent\/{\em :}
   \begin{enumerate}
   \item[(i)] $\slam$ is formally normal,
   \item[(ii)] there exists a sequence
$\{u_{n}\}_{n=-\infty}^{\infty} \subseteq V$ such that
   \begin{align*}
\text{$u_{n-1} = \pa {u_n}$ and $|\lambda_{u_{n-1}}|=
|\lambda_{u_{n}}|$ for all $n \in \zbb$,}
   \end{align*}
and $\lambda_v = 0$ for all $v \in V \setminus
{\{u_{n}\colon n\in\zbb\}}$,
   \item[(iii)] $\slam \in \ogr{\ell^2(V)}$ and $\slam$ is
normal.
   \end{enumerate}
   \end{pro}
   \begin{proof}
(i)$\Rightarrow$(iii) It follows from
\cite[Proposition 3.4.3]{j-j-s} that
   \begin{align} \label{fn1+}
\text{$\escr \subseteq \dz{\slam^*\slam}$ and $\slam^*
\slam e_u = \|\slam e_u\|^2 e_u$ for all $u \in V$.}
   \end{align}
Assertions (ii) and (iv) of Proposition \ref{bas}
imply that
   \begin{align}  \label{fn2+}
\text{$\escr \subseteq \dz{\slam \slam^*}$ and $\slam
\slam^* e_u \overset{\eqref{eu}\&\eqref{sl*}}= \sum_{v
\in \dzi{\pa{u}}} \lambda_v \overline{\lambda_u} e_v$
for all $u \in V^\circ$.}
   \end{align}
In view of \eqref{fn1+} and \eqref{fn2+}, we have
   \begin{align} \label{fn3+}
\escr \subseteq \xx:=\dz{\slam^* \slam} \cap \dz{\slam
\slam^*}.
   \end{align}
The formal normality of $\slam$ yields
$\is{\slam^*\slam f}{f} = \is{\slam\slam^* f}{f}$ for
all $f \in \xx$. Hence, in view of the polarization
formula, we have $\is{\slam^*\slam f}{g} =
\is{\slam\slam^* f}{g}$ for all $f, g \in \xx$. Since,
by \eqref{fn3+}, the vector space $\xx$ is dense in
$\ell^2(V)$, we obtain $\slam^* \slam e_u = \slam
\slam^* e_u$ for all $u \in V$. This, combined with
\eqref{fn1+} and \eqref{fn2+}, shows that
   \begin{align} \label{normal1}
\|\slam e_u\|^2 e_u = |\lambda_u|^2 e_u + \sum_{v \in
\dzi{\pa{u}} \setminus \{u\}} \lambda_v
\overline{\lambda_u} e_v, \quad u \in V^\circ.
   \end{align}

Note that if $u \in V^\circ$ is such that $\|\slam
e_u\|> 0$, then by \eqref{normal1} we have $\|\slam
e_u\|=|\lambda_u|$, $\lambda_u \neq 0$ and $\lambda_v
= 0$ for all $v \in \dzi{\pa{u}} \setminus \{u\}$,
which, in view of \eqref{eu}, implies that
   \begin{align*}
\|\slam e_{\paa(u)}\| = \|\slam e_u\| > 0.
   \end{align*}
Using an induction argument, we show that the
following implication holds for all $u \in V$ and $m
\in \zbb_+$ such that $\paa^{m}(u)\in V$:
   \begin{align}  \label{normal2}
\text{if $\|\slam e_u\|> 0$, then $\|\slam
e_{\paa^k(u)}\| = \|\slam e_u\|$ for all $k =0,\ldots,
m$.}
   \end{align}
Now we prove that if $u_1, u_2 \in V$ are such that
$\|\slam e_{u_1}\|> 0$ and $\|\slam e_{u_2}\|> 0$,
then $\|\slam e_{u_1}\| = \|\slam e_{u_2}\|$. Indeed,
by \cite[Proposition 2.1.4]{j-j-s}, there exists $u
\in V$ such that $u_1,u_2 \in \des{u}$. It follows
from \eqref{num4} and \eqref{num3} that there are
$m_1, m_2 \in \zbb_+$ such that
$\paa^{m_1}(u_1)=u=\paa^{m_2}(u_2)$. This fact
combined with \eqref{normal2} leads to $\|\slam
e_{u_1}\| = \|\slam e_{u}\| = \|\slam e_{u_2}\|$. This
implies that $\sup_{v\in V}\|\slam e_{v}\| < \infty$,
which together with assertions (i) and (iii) of
Proposition \ref{bas} yields $\slam \in
\ogr{\ell^2(V)}$. Hence $\slam$ is a bounded normal
operator.

(iii)$\Rightarrow$(i) Evident.

(ii)$\Rightarrow$(iii) It follows from (ii) that for
any $u\in V$ the set $\{v \in \dzi{u} \colon \lambda_v
\neq 0\}$ has at most one element. Consequently, again
by (ii), $\sup_{u \in V} \sum_{v \in \dzi{u}}
|\lambda_v|^2 = |\lambda_{u_1}|^2 < \infty$, which,
combined with Proposition \ref{bas}\,(iii), implies
that $\slam \in \ogr{\ell^2(V)}$. Now applying
\eqref{eu} and \eqref{sl*} separately to $u \in
{\{u_{n}\colon n\in\zbb\}}$ and $u \in V \setminus
{\{u_{n}\colon n\in\zbb\}}$, we verify that $\slam^*
\slam e_u = \slam \slam^* e_u$ for all $u \in V$,
which yields the normality of $\slam$.

(iii)$\Rightarrow$(ii) Apply \cite[Lemma 8.1.5]{j-j-s}
(this is the only case in which we use the assumption
that $\slam$ is a nonzero operator).
   \end{proof}
Combining Proposition \ref{normal} with
\cite[Proposition 8.1.6]{j-j-s}, we see that the only
directed tree admitting formally normal weighted
shifts with nonzero weights is isomorphic to $(\zbb,
\{(n,n+1)\colon n \in \zbb\})$.
   \section{\label{nevws}Modeling normal extensions
on weighted shifts} In this final section we will
discuss the following question:\ under what
circumstances can a normal extension of a subnormal
weighted shift on a directed tree $\tcal$ be modeled
as a weighted shift on a directed tree $\hat \tcal$
(no relationship between $\tcal$ and $\hat\tcal$ is
required). Here and in what follows, by the {\em
unilateral shift} on $\ell^2(\zbb_+)$ (respectively:\
the {\em bilateral shift} on $\ell^2(\zbb)$) we mean
the weighted shift on the directed tree $(\zbb_+,
\{(n,n+1)\colon n \in \zbb_+\})$ (respectively:\
$(\zbb, \{(n,n+1)\colon n \in \zbb\}$)) with all
weights equal to $1$. These two particular directed
trees are denoted simply by $\zbb_+$ and $\zbb$,
respectively.
   \begin{lem} \label{u+0}
Let $\slam$ be a nonzero subnormal weighted shift on a
directed tree $\tcal$. Suppose $\slam$ has a normal
extension $N$ which is a weighted shift on a directed
tree $\hat{\tcal}=(\hat V, \hat E)$ $($we do not
assume that $\tcal$ is a directed subtree of
$\hat\tcal$$)$. Then $N \in \ogr{\ell^2(\hat V)}$ and
$N = \alpha U \oplus 0$, where $\alpha$ is a positive
real number, $U$ is a unitary operator which is
unitarily equivalent to the bilateral shift on
$\ell^2(\zbb)$ and $0$ is the zero operator on
$\ell^2(\hat V) \ominus \dz{U}$. Moreover, $\ob{\slam}
\subseteq \dz{U}$.
   \end{lem}
   \begin{proof}
Denote by $\{\hat \lambda_v\}_{v \in \hat V^\circ}$
the weights of $N$. It follows from Proposition
\ref{normal} that $N$ is a bounded operator on
$\ell^2(\hat V)$ and that there exist a positive real
number $\alpha$ and a sequence
$\{u_{n}\}_{n=-\infty}^{\infty} \subseteq \hat V$ such
that
   \begin{align} \label{pr2-1}
& \text{$u_{n-1} = \paa_{\hat \tcal}(u_n)$ for all $n
\in \zbb$,}
   \\
& \text{$\hat \lambda_v = 0$ for all $v \in \hat V
\setminus {\{u_{n}\colon n\in\zbb\}}$,} \label{pr2-2}
   \\
& \text{$|\hat\lambda_{u_{n}}|= \alpha$ for all $n \in
\zbb$,} \label{pr2-3}
   \end{align}
where $\paa_{\hat \tcal}(\cdot)$ refers to the
directed tree $\hat \tcal$. Set $X=\{u_n\colon n \in
\zbb\}$ and $Y=\hat V \setminus X$. We deduce from
\eqref{eu}, \eqref{pr2-1} and \eqref{pr2-2} that the
spaces $\ell^2(X)$ and $\ell^2(Y)$ are invariant for
$N$ (and thus $N=N|_{\ell^2(X)} \oplus
N|_{\ell^2(Y)}$), $N|_{\ell^2(Y)}=0$ and $N\hat
e_{u_n}= \hat \lambda_{u_{n+1}} \hat e_{u_{n+1}}$ for
all $n \in \zbb$, where $\{\hat e_{v}\}_{v \in \hat
V}$ is the canonical orthonormal basis of $\ell^2(\hat
V)$. Applying \eqref{pr2-3} and \cite[Corollary 1, p.\
52]{shi} we get the required decomposition $N = \alpha
U \oplus 0$. Hence, we have $\ob{\slam}\subseteq
\ob{N}=\ob{U}=\dz{U}$, which completes the proof.
   \end{proof}
Regarding Lemma \ref{u+0}, we note that bounded
subnormal operators with normal extensions of the form
$U \oplus 0$, where $U$ is a unitary operator, have
been characterized in \cite{c-st}.

Now we show that the only nonzero subnormal weighted
shifts on directed trees with nonzero weights whose
normal extensions can be modeled as weighted shifts on
directed trees are those that are unitarily equivalent
to either positive scalar multiples of the bilateral
shift on $\ell^2(\zbb)$ or positive scalar multiples
of ``small'' perturbations of the unilateral shift on
$\ell^2(\zbb_+)$.
   \begin{thm}\label{normwsno}
Let $\slam$ be a nonzero subnormal weighted shift on a
directed tree $\tcal$ with nonzero weights
$\lambdab=\{\lambda_v\}_{v \in V^\circ}$. Suppose
$\slam$ has a normal extension $N$ which is a weighted
shift on a directed tree $\hat{\tcal}=(\hat V, \hat
E)$ $($we do not assume that $\tcal$ is a directed
subtree of $\hat\tcal$$)$. Then the directed tree
$\tcal$ is isomorphic to either $\zbb$ or $\zbb_+$. In
the former case, $\slam$ is unitarily equivalent to a
positive scalar multiple of the bilateral shift on
$\ell^2(\zbb)$. In the latter case, $\slam$ is
unitarily equivalent to a positive scalar multiple of
a unilateral weighted shift on $\ell^2(\zbb_+)$ with
weights $\{\vartheta, 1, 1, 1, \ldots\}$, where
$\vartheta \in (0,1]$.
   \end{thm}
   \begin{proof}
In view of Lemma \ref{u+0}, the operator $N$ (and
consequently $\slam$) is bounded and $N = \alpha U
\oplus 0$, where $\alpha$ and $U$ are as in Lemma
\ref{u+0}. Without loss of generality, we can assume
that $\alpha=1$. It follows from Lemma \ref{u+0} that
there exists a unitary isomorphism $W\colon \dz{U} \to
L^2(\tbb)$ such that $WU=MW$, where $L^2(\tbb)$ is the
Hilbert space of all square summable Borel functions
on $\tbb:=\{z\in \cbb\colon |z|=1\}$ with respect to
the normalized Lebesgue measure $m$ on $\tbb$ that is
given by
   \begin{align*}
m(\sigma)=\frac {1}{2\pi} \int_0^{2\pi} \chi_{\sigma}
(\E^{\I t}) \D t, \quad \sigma \text{ - Borel subset
of } \tbb,
   \end{align*}
($\chi_\sigma$ stands for the characteristic function
of $\sigma$), and $M\in \ogr{L^2(\tbb)}$ is defined by
$(Mf)(z)=zf(z)$ a.e.\ $[m]$ for $f \in L^2(\tbb)$.

By \cite[Proposition 5.1.1]{j-j-s}, the directed tree
$\tcal$ is leafless. To prove that the directed tree
$\tcal$ is isomorphic to either $\zbb$ or $\zbb_+$, it
is enough to show that $\tcal$ has no branching
vertex. Suppose that, contrary to our claim, $\tcal$
has a branching vertex $u \in V$. Then, by
\eqref{num3+}, there exist $u_1,u_2 \in V$ such that
   \begin{align} \label{desdes0}
\des{u_1} \cap \des{u_2} = \varnothing.
   \end{align}
It follows from \cite[Lemma 6.1.1]{j-j-s} that
$\slam^n e_{u_j} \in \ell^2(\dzin{n}{u_j})$ for all $n
\in \zbb_+$ and $j=1,2$. This fact, combined with
\eqref{num3} and \eqref{desdes0}, implies that
   \begin{align} \label{ort3}
\text{the vectors $\{\slam^k e_{u_1}\colon k\in
\zbb_+\} \cup \{\slam^l e_{u_2}\colon l\in \zbb_+\}$
are pairwise orthogonal.}
   \end{align}
Since $\ob{\slam} \subseteq \dz{U}$, we see that
$\slam e_{u_j} \in \dz{U}$ for $j=1,2$. Set $f_j =
W\slam e_{u_j} \in L^2(\tbb)$ for $j=1,2$. In view of
Proposition \ref{bas}\,(v), the operator $\slam$ is
injective. Hence $\|f_j\|>0$ for $j=1,2$. It follows
from the equality $N = U \oplus 0$ that
   \begin{align}  \label{ort4}
M^n f_j = W U^n(\slam e_{u_j}) = W \slam^{n+1}
e_{u_j}, \quad n \in \zbb_+, \, j=1,2.
   \end{align}
Combining conditions \eqref{ort3} and \eqref{ort4}, we
deduce that the vectors $\{M^k f_1\colon k\in \zbb_+\}
\cup \{M^l f_2\colon l\in \zbb_+\}$ are pairwise
orthogonal. Thus the following equalities hold for all
$k \in \nbb$, $l \in \zbb_+$ and $j=1,2$:
   \begin{gather*}
\int_{\tbb} z^k |f_j(z)|^2 \D m (z) = \is{M^k
f_j}{f_j} = 0 = \is{f_j}{M^k f_j} = \int_{\tbb} \bar
z^k |f_j(z)|^2 \D m (z),
   \\
\int_{\tbb} z^l f_1(z) \overline{f_2(z)} \D m (z) =
\is{M^l f_1}{f_2} = 0 = \is{f_1}{M^l f_2} =
\int_{\tbb} \bar z^l f_1(z) \overline{f_2(z)} \D m
(z).
   \end{gather*}
This yields
   \begin{gather}  \label{ort1}
\int_{\tbb} z^k |f_j(z)|^2 \D m (z) = 0, \quad k\in
\zbb \setminus \{0\}, \, j=1,2,
   \\   \label{ort2}
\int_{\tbb} z^l f_1(z) \overline{f_2(z)} \D m (z) = 0,
\quad l \in \zbb.
   \end{gather}
It follows from \eqref{ort1} that for every complex
trigonometric polynomial $p$ on $\tbb$,
   \begin{align*}
\int_{\tbb} p \, |f_j|^2 \D m = \|f_j\|^2 \int_{\tbb}
p \D m, \quad j=1,2.
   \end{align*}
Since complex trigonometric polynomials are uniformly
dense in the Banach space of continuous functions on
$\tbb$, we infer from the Riesz representation theorem
(cf.\ \cite[Theorem 6.19]{Rud}) that $\int_{\sigma}
|f_j|^2 \D m = \|f_j\|^2 m(\sigma)$ for all Borel
subsets $\sigma$ of $\tbb$. This implies that
   \begin{align} \label{ny1}
|f_j| = \|f_j\| \quad \text{a.e.\ $[m]$.}
   \end{align}
A similar argument applied to \eqref{ort2} yields
   \begin{align}  \label{ny2}
f_1 \overline{f_2} = 0 \quad \text{a.e.\ $[m]$.}
   \end{align}
Combining \eqref{ny1} and \eqref{ny2} with the fact
that $\|f_j\|>0$ for $j=1,2$, we conclude that
$m(\tbb)=0$, which contradicts the fact that
$m(\tbb)=1$. This shows that the directed tree $\tcal$
is isomorphic to either $\zbb$ or $\zbb_+$.

First we consider the case in which $\tcal$ is
isomorphic to $\zbb$. Without loss of generality, we
can assume that $\tcal$ coincides with the directed
tree $\zbb$. Since the weights of $\slam$ are nonzero,
we get $\mathscr{E}_{\zbb} \subseteq \ob{\slam}
\subseteq \dz{U}$, which yields $\ell^2(\zbb)
\subseteq \dz{U}$. Hence $U$ is a unitary extension of
$\slam$. This implies that $\slam $ is an isometric
bilateral weighted shift on $\ell^2(\zbb)$. As a
consequence of \cite[Corollary 1, p.\ 52]{shi}, the
operator $\slam$ is unitarily equivalent to the
bilateral shift on $\ell^2(\zbb)$.

Consider now the case in which $\tcal$ is isomorphic
to $\zbb_+$. Again without loss of generality, we can
assume that $\tcal$ coincides with the directed tree
$\zbb_+$. Since the weights of $\slam$ are nonzero, we
obtain $\lin \{e_n \colon n \in \nbb\} \subseteq
\ob{\slam} \subseteq \ell^2(\nbb)$, and so
$\ell^2(\nbb) = \overline{\ob{\slam}} \subseteq
\dz{U}$. As $U$ is a unitary operator and $U \oplus 0$
extends $\slam$, we deduce that
$\slam|_{\overline{\ob{\slam}}}$ is an isometry.
Hence, we have
   \begin{align*}
 \|\slam^n e_{1}\|^2 =
 \frac{1}{|\lambda_1|^2}\|\slam^n \slam e_{0}\|^2 =
 \frac{1}{|\lambda_1|^2} \|\slam e_{0}\|^2 = 1, \quad
 n \in \zbb_+.
   \end{align*}
This implies that $\{\|\slam^n e_1\|^2\}_{n=0}^\infty$
is a Stieltjes moment sequence with a representing
measure $\delta_1$ ($\delta_1$ is the Borel
probability measure on $[0,\infty)$ concentrated at
the point $1$). Since $\slam$ is subnormal, we deduce
that $\{\|\slam^{n} e_0\|^2\}_{n=0}^\infty$ is a
Stieltjes moment sequence (see e.g., \cite[Theorem
6.1.3]{j-j-s}). Hence, by applying \cite[Lemma
6.1.10]{j-j-s} to $u=0$, we get $|\lambda_1| \Le 1$.
The fact that $\slam|_{\overline{\ob{\slam}}}$ is an
isometry yields $|\lambda_n|=1$ for all $n \Ge 2$.
Using \cite[Corollary 1, p.\ 52]{shi}, we conclude
that $\slam$ is unitarily equivalent to a unilateral
weighted shift on $\ell^2(\zbb_+)$ with weights
$\{\vartheta, 1, 1, 1, \ldots\}$, where $\vartheta \in
(0,1]$. This completes the proof.
   \end{proof}
   \begin{rem}  \label{rem}
We show how to model a normal extension of a nonzero
subnormal weighted shift $\slam$ on a directed tree
$\tcal$ with nonzero weights by means of a weighted
shift on some directed tree. As in the proof of
Theorem \ref{normwsno} we consider only the case of
$\alpha=1$. By this theorem we have only two
possibilities:\ either the directed tree $\tcal$ is
isomorphic to $\zbb$ and $\slam$ is normal (and so
$\slam$ is the required model), or the directed tree
$\tcal$ is isomorphic to $\zbb_+$ and $\slam$ is a
unilateral weighted shift on $\ell^2(\zbb_+)$ with
weights $\{\vartheta, 1, 1, 1, \ldots\}$, where
$\vartheta \in (0,1]$. In the latter case, we fix
$\omega \notin \zbb$ and define the directed tree
$\hat\tcal = (\hat V, \hat E)$ by
   \begin{align*}
\hat V & = \{\omega\} \cup \zbb \quad \text{and} \quad
\hat E = \{(0,\omega)\} \cup \{(n,n+1)\colon n \in
\zbb\}.
   \end{align*}
Then $\hat \tcal$ is rootless and $0$ is a unique
branching vertex of $\hat \tcal$. Let $N$ be the
weighted shift on $\hat \tcal$ with weights $\{\hat
\lambda_v\}_{v\in \hat V}$ given by $\hat \lambda_v =
0$ for $v = \omega$ and $\hat \lambda_v = 1$ for $v
\in \zbb$. By Proposition \ref{normal}, $N$ is a
bounded normal operator on $\ell^2(\hat V)$. Define
the sequence $\{\tilde e_n\}_{n=0}^\infty$ in
$\ell^2(\hat V)$ by
   \begin{align*}
\tilde e_n =
   \begin{cases}
\sqrt{(1-\vartheta^2)} \, \hat e_{\omega} + \vartheta
\hat e_0 & \text{ for } n=0,
   \\[.5ex]
\hat e_n & \text{ for } n \in \nbb.
   \end{cases}
   \end{align*}
Denote by $\hh$ the closure of $\lin\{\tilde
e_n\}_{n=0}^\infty$ in $\ell^2(\hat V)$. Then
$\{\tilde e_n\}_{n=0}^\infty$ is an orthonormal basis
of $\hh$. It is clear that $N \tilde e_0 = \vartheta
\tilde e_1$ and $N \tilde e_n = \tilde e_{n+1}$ for
all $n \in \nbb$. As a consequence, $N(\hh) \subseteq
\hh$ and the operator $N|_{\hh}$ is unitarily
equivalent to a unilateral weighted shift on
$\ell^2(\zbb_+)$ with weights $\{\vartheta, 1, 1, 1,
\ldots\}$. Hence, $N$ is the required model of a
normal extension of $\slam$. It is worth noting that
the directed tree $\tcal$ is isomorphic to many
directed subtrees of $\hat\tcal$. However, if
$\vartheta \in (0,1)$, then $\tcal$ could not be
regarded as a directed subtree of $\hat \tcal$.
Indeed, otherwise $V \subseteq \zbb \subseteq \hat V$
and so $\slam$ is isometric as the restriction to
$\ell^2(V)$ of the unitary operator
$N|_{\ell^2(\zbb)}$. This contradicts the fact that
for $\vartheta \in (0,1)$, $\slam$ is not an isometry.
   \end{rem}
   \subsection*{Acknowledgement} The final version of
this paper was written while the first and the third
authors visited Kyungpook National University during
the spring of 2011. They wish to thank the faculty and
the administration of this unit for their warm
hospitality.
   \bibliographystyle{amsalpha}
   
   \end{document}